\documentclass[12pt,twoside]{amsart}
\usepackage{amssymb}
\usepackage{amscd}

\title[Multiplication and vanishing]
{Multiplication maps and vanishing theorems for 
toric varieties}
\author{Osamu Fujino} 
\date{2006/11/30}
\address{Graduate School of Mathematics\\ 
 Nagoya University, Chikusa-ku Nagoya 464-8602 Japan}
\email{fujino@math.nagoya-u.ac.jp}
\newcommand{\Hom}[0]{{\operatorname{Hom}}}
\newcommand{\codim}[0]{{\operatorname{codim}}}
\newtheorem{thm}{Theorem}[section]
\newtheorem{lem}[thm]{Lemma}
\newtheorem{cor}[thm]{Corollary}
\newtheorem{prop}[thm]{Proposition}

\theoremstyle{definition}
\newtheorem{ex}[thm]{Example}
\newtheorem{defn}[thm]{Definition}
\newtheorem{rem}[thm]{Remark}
\newtheorem*{ack}{Acknowledgments}       
\newtheorem{say}[thm]{}
\begin{document}
\bibliographystyle{amsalpha+}

\begin{abstract} 
We use multiplication maps to give a characteristic-free approach 
to vanishing theorems on toric varieties. 
Our approach is very elementary but is enough powerful 
to prove vanishing theorems. 
\end{abstract}
\maketitle
\tableofcontents

\section{Introduction}\label{sec1} 
The main purpose of this paper is to understand 
various vanishing theorems on toric varieties 
through multiplication maps. 
We give an elementary and unified approach to vanishing theorems on 
toric varieties.  
The following theorem is the main theorem of this paper. 
Some important special cases were already investigated in 
various papers. See, for example, \cite[7.5.2.~Theorem]{danilov}, 
\cite[Section 7]{bc}, 
\cite[Theorem 5]{btlm}, and \cite[Section 2]{mustata}. 

\begin{thm}[Main theorem I]\label{main} 
Let $X$ be a toric variety defined over a field $k$ of 
arbitrary characteristic and $B$ 
a reduced torus invariant Weil divisor on $X$. 
Let $L$ be a line bundle on $X$. If 
$H^i(X, \widetilde {\Omega}^a_{X}(\log B)\otimes L^l)=0$ for 
some positive integer $l$, 
then $H^i(X, \widetilde {\Omega}^a_{X}(\log B)\otimes L)=0$. 
In particular, if $X$ is projective and 
$L$ is ample, then $H^i(X, \widetilde {\Omega}^a_{X}(\log B)\otimes L)=0$ 
for any $i>0$. 
\end{thm} 

Before we go further, let us recall the definition of $\widetilde \Omega^a_X(\log B)$ 
(cf.~\cite[15.2]{danilov}).  

\begin{defn}\label{12}
Let $W$ be any Zariski open set of $X$ such that $W$ is non-singular and 
$\codim _X(X\setminus W)\geq 2$. 
Furthermore, we assume that $B$ is a simple normal crossing divisor on 
$W$. On this assumption, 
$\Omega^a_{W}(\log B)$ 
is a well-defined locally free sheaf on $W$. 
Let $\iota:W\hookrightarrow X$ be the natural 
open immersion. 
Then we put $\widetilde \Omega^a_X(\log B)=\iota_*\Omega^a_W(\log B)$ for any 
$a\geq 0$. 
It is easy to see that the reflexive sheaf 
$\widetilde \Omega^a_{X}(\log B)$ on $X$ does not 
depend on the choice of $W$. 
Note that $B$ is a simple normal crossing 
divisor on $W$ if $W$ is a non-singular 
toric variety. 
If $B=0$, then we write $\widetilde \Omega^a_{X}=
\widetilde \Omega^a_X(\log B)$ for any $a\geq 0$. 
\end{defn}

The above theorem contains the 
following important vanishing theorems. 
If $B=0$, then we obtain the famous 
Bott type vanishing theorem for 
toric varieties. It was first claimed 
in \cite[7.5.2~Theorem]{danilov} without 
proof. See \cite[Theorem 5]{btlm}. 
The readers can find that this famous 
vanishing theorem is stated in the standard reference 
\cite[p.130]{odasensei} without proof. 

\begin{cor}[Bott, Danilov, $\cdots$]\label{bott}
Let $X$ be a projective toric variety over $k$ 
and $L$ an ample line bundle on $X$. 
Then $H^i(X, \widetilde \Omega^a_X\otimes L)=0$ for any $i>0$ and $a\geq 0$. 
\end{cor}
 
In the main theorem:~Theorem \ref{main}, if we put $a=\dim X$, then we obtain 
the toric version of Norimatsu type vanishing theorem. 
It is nothing but Musta\c t\u a's vanishing theorem in \cite[Corollary 2.5 (iii)]{mustata}. 
The readers can find the original formulation in 
Corollary \ref{honyaku} below. One of my motivations 
is to give an elementary proof to Musta\c t\u a's vanishing theorem. 

\begin{cor}[Norimatsu, Musta\c t\u a, $\cdots$]\label{nori}
Let $X$ be a projective toric variety over $k$ and $B$ a reduced torus invariant 
Weil divisor on $X$. Let $L$ be an ample line bundle on $X$. 
Then $H^i(X, \mathcal O_X(K_X+B)\otimes L)=0$ for any $i>0$. 
\end{cor}

Note that $K_X$ is the canonical divisor of $X$. 
It is well known that $\mathcal O_X(K_X)\simeq \mathcal O_X(-\sum _iD_i)$, 
where the summation $\sum _iD_i$ runs over 
all the torus invariant prime divisors on $X$. 
The final one is the Kodaira type vanishing theorem for toric varieties. 
It is sufficiently powerful in the toric geometry (see Section \ref{sec3}). 

\begin{cor}[Kodaira, $\cdots$]\label{kodaira}
Let $X$ be a projective toric variety over $k$ and $L$ an ample line bundle on $X$. 
Then $H^i(X, \mathcal O_X(K_X)\otimes L)=0$ for any $i>0$. 
\end{cor}

The next theorem is another main theorem of this paper. 
It contains the Kawamata-Viehweg type vanishing theorem obtained by Musta\c t\u a 
(see \cite[Corollary 2.5 (i) and (ii)]{mustata}). 
Our formulation is very similar to Musta\c t\u a's 
theorem:~\cite[Theorem 0.1]{mustata}, but is slightly different. 
We will 
quickly see the relationship between 
Musta\c t\u a's original statement and 
Theorem \ref{main2} in \ref{214}. 
His statement is a special case of our theorem (see Corollary 
\ref{musumusu}). 
See also Proposition \ref{vari}, 
where we will treat a variant of Theorem \ref{main2}.  

\begin{thm}[Main theorem II]\label{main2}
Let $X$ be a toric variety defined 
over a field 
$k$ of arbitrary characteristic and $D$ a torus invariant $\mathbb Q$-Weil 
divisor on $X$. Assume that $lD$ is an integral 
Weil divisor for some positive integer $l$. 
If $H^i(X, \mathcal O_X(lD))=0$ $($resp.~$H^i(X, \mathcal O_X(K_X+lD))=0$$)$, 
then we have 
$H^i(X, \mathcal O_X(\llcorner D\lrcorner ))=0$ $($resp.~$H^i(X, \mathcal 
O_X(K_X+\ulcorner D\urcorner ))=0$$)$. 
\end{thm} 

The following corollary easily follows from Theorem 
\ref{main2}. 
However, \cite[Theorem 0.1]{mustata} produced 
it only when $D$ is an ample $\mathbb Q$-Cartier 
divisor (see \cite[Corollary 2.5 (i) 
and (ii)]{mustata}). 

\begin{cor}[Kawamata-Viehweg, Musta\c t\u a, $\cdots$]\label{kavi} 
Let $X$ be a complete toric variety over $k$ and $D$ a nef 
$\mathbb Q$-Cartier torus invariant $\mathbb Q$-Weil divisor on $X$ 
with 
the Iitaka dimension 
$\kappa (X, D)=\kappa$. 
Then we obtain $H^i(X, \mathcal O_X(\llcorner D\lrcorner))=0$ for 
$i\ne 0$ and 
$H^i(X, \mathcal O_X(K_X+\ulcorner D\urcorner ))=0$ for 
$i\ne n-\kappa$, where $n=\dim X$. 
\end{cor} 

We note that for a $\mathbb Q$-Weil divisor 
$D=\sum _{j=1}^{r}d_jD_j$ on $X$, we 
define the round-up 
$\ulcorner D\urcorner =\sum _{j=1}^{r}\ulcorner d_j\urcorner 
D_j$ (resp.~the round-down 
$\llcorner D\lrcorner =\sum _{j=1}^{r}\llcorner d_j\lrcorner D_j$), 
where for any real number $x$, $\ulcorner x\urcorner $ (resp.~$\llcorner 
x\lrcorner$) is the integer defined by $x\leq \ulcorner x\urcorner <x+1$ (resp.~$x-1
<\llcorner x\lrcorner \leq x$). 
The fractional part $\{D\}$ of the 
$\mathbb Q$-Weil divisor $D$ denotes $D-\llcorner D\lrcorner$. 

We summarize the contents of this 
paper:~In Section \ref{sec2}, we will prove Theorem \ref{main} 
and Theorem \ref{main2}. 
The main ingredient of our proof is the 
{\em{multiplication map}}. 
It is a mystery that no standard references on the 
toric geometry treat the multiplication map systematically. 
Let us introduce the $l$ times multiplication 
map for toric varieties. 
We consider $\mathbb P^n$ and 
a finite surjective 
morphism $F:\mathbb P^n\to \mathbb P^n: [X_0:\cdots: X_n]\mapsto 
[X^l_0:\cdots:X^l_n]$. 
It is the simplest example of $l$ times 
multiplication maps for projective 
toric varieties. 
On the big torus $T\subset \mathbb P^n$, the 
restriction $F_T:=F|_T:T\to T$ is nothing 
but the group homomorphism expressed 
by $(x_1, \cdots, x_n)\mapsto 
(x_1^l, \cdots, x_n^l)$. For 
an arbitrary $n$-dimensional 
toric variety $X$, $F_T:T\to T$ naturally 
extends to a finite surjective toric morphism 
$F:X\to X$. 
We call this $F:X\to X$ the {\em{$l$ times 
multiplication map}} of $X$. 
I believe that the multiplication map 
will play important roles in the toric geometry. 
Here, I will show its usefulness by proving 
various vanishing theorems. 
Our approach is very 
elementary but is sufficiently powerful to prove 
vanishing theorems. 
Related topics are in Section 7 of \cite{arapura}. 
We do not use Frobenius morphisms (cf.~\cite{btlm} and \cite{blickle}) 
nor Cox's homogeneous coordinate rings (cf.~\cite{mustata}). 
We do not need any cumbersome combinatorial arguments 
nor the Hodge theory (cf.~\cite{bc}). 
We recommend the readers to compare our proof with the 
others (cf.~\cite{bc}, \cite{btlm}, \cite{mustata}, etc.). 
In Section \ref{varivari}, we consider slight generalizations 
of the main theorems:~Theorem \ref{main} and Theorem 
\ref{main2}. Let $\mathcal E$ be a reflexive sheaf on $X$. 
Roughly speaking, we treat 
$(\mathcal E\otimes {\widetilde \Omega}^a(\log B))^{**}$, 
$(\mathcal E\otimes \mathcal O_X(\llcorner 
D\lrcorner))^{**}$, 
and $(\mathcal E\otimes \mathcal O_X(K_X+\ulcorner 
D\urcorner))^{**}$ 
instead of ${\widetilde \Omega}^a(\log B)$, 
$\mathcal O_X(\llcorner D\lrcorner)$, 
and $\mathcal O_X(K_X+\ulcorner D\urcorner)$ 
respectively. We note 
that $\mathcal E$ is not assumed to be equivariant. 
In Section \ref{sec3}, we will treat Koll\'ar's injectivity 
theorem for toric varieties. For toric varieties, it easily follows from 
the Kodaira type vanishing theorem. 
In Section \ref{sec4}, which is an appendix, 
we will state relative vanishing 
theorems explicitly for future uses. 

We note that our reference list does not cover all the papers 
treating the related topics. We apologize in advance to the colleagues whose 
works are not appropriately mentioned in this paper. 

\begin{ack} 
I would like to express my gratitude to Professors Kazuhiro Fujiwara and 
Takeshi Abe for giving me much advice and encouraging me 
during the preparation of this paper. I thank 
Doctor Hiroshi Sato for valuable 
discussions and Professor Noboru Nakayam for 
pointing out a little mistake. I was 
partially supported by The Sumitomo 
Foundation and by the Grant-in-Aid for 
Young Scientists (A) $\sharp$17684001 
from JSPS. 
I would like to thank Professor Donu Arapura and
Doctor Sam Payne, who gave me comments by e-mails
after I circulated the first version of this paper.
Professor Donu Arapura informed me of his recent paper \cite{arapura}.
Doctor Sam Payne sent me his private notes on vanishing theorems
for toric varieties. 
The discussions with him helped me revise this paper. 
\end{ack}

Let $k$ be a fixed field of arbitrary characteristic $p$ ($p$ may be zero). 
In this paper, everything is defined over $k$. 
We do not assume that $k$ is algebraically closed. 

\section{Multiplication maps and 
vanishing theorems}\label{sec2}
We fix our notation and define the multiplication map. 

\begin{say}\label{21}
Let $N\simeq \mathbb Z^n$ be a lattice 
and $M=\Hom _{\mathbb Z}(N, \mathbb Z)$ the dual lattice. 
For a fan $\Delta$ in $N_{\mathbb R}=N\otimes _{\mathbb Z}\mathbb R$, we 
have the associated toric variety $X=X(\Delta)$. 
We put $N'=\frac{1}{l}N$ and $M'=\Hom _{\mathbb Z}(N', \mathbb Z)$ for 
any positive integer $l$. 
We note that $M'=lM$. 
Since $N_{\mathbb R}=N'_{\mathbb R}$, $\Delta$ is also a fan in 
$N'_{\mathbb R}$. 
We write $\Delta'$ to express the fan $\Delta$ in $N'_{\mathbb R}$. 
Let $X'=X(\Delta')$ be the associated toric variety. 
We note that $X\simeq X'$ as toric varieties. 
We consider 
the natural inclusion $\varphi:N\to N'$. 
Then $\varphi$ induces a finite surjective toric morphism $F:X\to X'$. 
We call it the {\em{$l$ times multiplication map of $X$}}. 
The following is the most important example of $l$ 
times multiplication maps. 
\begin{ex}
The finite surjective morphism 
$F:\mathbb A^n\to \mathbb A^n$ given by $(a_1, \cdots, a_n)\mapsto 
(a_1^l, \cdots, a_n^l)$ is the $l$ times multiplication 
map of $\mathbb A^n$. 
\end{ex}
\end{say} 

Let us start the proof of the main theorem I:~Theorem \ref{main}. 

\begin{say}\label{22}
Let $\mathcal A$ be an object on $X$. Then we write $\mathcal A'$ 
to indicate the corresponding 
object on $X'$. 
Let $T$ be the big torus of $X$. 
We construct a split injection $\Omega^1_{T'}\to F_*\Omega^1_{T}$. 
Note that $\Omega ^1_T$ is nothing but a $k[M]$-module 
$M\otimes _\mathbb Z k[M]$. 
\end{say}

We recall the toric description of $\Omega^1_T$ more precisely. 
For the details, see \cite[\S 4]{danilov} and \cite{ishida}. 

\begin{say}\label{23} 
By choosing a base suitably, we have 
$k[M]\simeq k[x_1, x^{-1}_1, \cdots, x_n, x^{-1}_n]$. 
We can write $x^m=x^{m_1}_1 x^{m_2}_2\cdots x^{m_n}_n$ for 
$m=(m_1, \cdots, m_n)\in \mathbb Z^n=M$. 
Then we have the isomorphism of $k[M]$-modules 
$M\otimes _{\mathbb Z}k[M]\to H^0(T, \Omega^1_T)$ induced by 
$m\otimes 
x^{\widetilde m}\mapsto \frac{dx^m}{x^m}\cdot 
x^{\widetilde m}=x^{\widetilde m-m}dx^m$, 
where $m, \widetilde m\in \mathbb Z^n=M$. 
Note that $\wedge ^a M\otimes _{\mathbb Z}k[M]\simeq 
H^0(T, \Omega^a_T)$ as $k[M]$-modules for any $a\geq 0$. 
\end{say} 

We go back to the proof of the main theorem. 

\begin{say}\label{24}
Therefore, $F_*\Omega^1_T$ corresponds to a $k[M']$-module $M\otimes 
_{\mathbb Z}k[M]$. 
We consider the $k[M']$-module homomorphism 
$M'\otimes _{\mathbb Z}k[M']\to M\otimes _{\mathbb Z}k[M]$ induced 
by $m_{\alpha}\otimes x^{m_{\beta}}\mapsto m_{\alpha}\otimes x^{lm_{\beta}}$. 
This gives an injection $\Omega^1_{T'}\to F_*\Omega^1_T$. 
We also consider the $k[M']$-module homomorphism 
$M\otimes _{\mathbb Z}k[M]\to M'\otimes _{\mathbb Z}k[M']$ obtained 
from $m_\alpha\otimes x^{m_\gamma}\mapsto m_\alpha\otimes x^{m_\beta}$ if 
$m_\gamma=lm_\beta$ and $m_\alpha\otimes x^{m_\gamma}\mapsto 0$ otherwise. 
By this homomorphism, the above injection 
$\Omega^1_{T'}\to F_*\Omega^1_T$ 
splits. 
We can generalize the above construction to $\wedge ^{a}M'\otimes 
_{\mathbb Z}k[M']$ and $\wedge ^{a}M\otimes _{\mathbb Z}k[M]$. 
More precisely, we consider the $k[M']$-module 
homomorphisms $\wedge ^aM'\otimes _{\mathbb Z}k[M']
\to \wedge ^aM\otimes _{\mathbb Z}k[M]$ given by 
$m_{\alpha_1}\wedge \cdots \wedge 
m_{\alpha_a}\otimes x^{m_{\beta}}\mapsto 
m_{\alpha_1}\wedge \cdots \wedge 
m_{\alpha_a}\otimes x^{lm_{\beta}}$, 
and $\wedge ^a M\otimes _{\mathbb Z}k[M]\to 
\wedge^aM'\otimes _{\mathbb Z}k[M']$ induced by 
$m_{\alpha_1}\wedge \cdots \wedge 
m_{\alpha_a}\otimes x^{m_{\gamma}}\mapsto 
m_{\alpha_1}\wedge \cdots \wedge 
m_{\alpha_a}\otimes x^{m_{\beta}}$ if $m_{\gamma}=lm_{\beta}$ 
and $m_{\alpha_1}\wedge \cdots \wedge 
m_{\alpha_a}\otimes x^{m_{\gamma}}\mapsto 0$ otherwise. 
So, 
we obtain split injections 
$\Omega^a_{T'}\to F_*\Omega^a_T$ for any $a\geq 0$. 
\end{say}
\begin{say}\label{26}
Let $\Delta_V$ be the fan in $N_{\mathbb R}$ that is 
obtained from $\Delta$ by removing the cones with 
dimensions $\geq 2$. 
Then $V=X(\Delta_V)$ is a non-singular toric variety such that 
$\codim _X(X\setminus V)\geq 2$. 
Let $B$ be a reduced torus invariant Weil divisor on $X$. 
Then we can construct split injections 
$\psi: \Omega^a_{V'}(\log B')\to F_*\Omega^a_V(\log B)$ 
for all $a\geq 0$, which are induced by 
$\Omega^a_{T'}\to F_*\Omega^a_T$. 
Note that we can see $\Omega^a_W(\log B)\subset \Omega^a_T$ 
for each affine toric open set $W$ of $V$. 
To check that $\psi: \Omega^a_{V'}(\log B')\to F_*\Omega^a_V(\log B)$ 
is a split injection, it is sufficient to check 
it on $U=k\times (k^{\times})^{n-1}\subset V$ since 
$V$ is covered by finitely many 
$k\times (k^{\times})^{n-1}$. On the open set 
$U$, it is easy to see that 
$\psi$ is a split injection by direct local computations. 
\end{say}
\begin{say}\label{27}
Let $\iota: V\hookrightarrow X$ be the natural open immersion. 
Since the following diagram 
$$
\begin{CD}
V@>{\iota}>>X\\ 
@V{F}VV @VV{F}V \\ 
V'@>{\iota'}>> X'
\end{CD}
$$ 
is commutative, 
we obtain split injections 
$\widetilde \psi =\iota'_*\psi:
\widetilde \Omega^a_{X'}(\log B')\to F_*\widetilde \Omega^a_{X}(\log B)$ 
for all $a$. 
Note that 
$\widetilde \Omega^a_X(\log B)=\iota_*\Omega^a_{V}(\log B)$ by Definition 
\ref{12}. 
\end{say}
\begin{say}\label{28}
Let $L$ be a line bundle on $X$. Since $L\simeq \mathcal O_X(G)$ 
for some torus invariant Cartier divisor $G$, we can see that 
$F^*L'\simeq L^l$. 
By combining these results, 
\begin{align*}
H^i(X, \widetilde \Omega^a_X(\log B)\otimes L)&\simeq 
H^i(X', \widetilde \Omega^a_{X'}(\log B')\otimes L')
\\ 
&\subset 
H^i(X', F_*\widetilde \Omega^a_X(\log B)\otimes L')\\&\simeq 
H^i(X, \widetilde \Omega^a_X(\log B)\otimes L^l). 
\end{align*} 
This inclusion and Serre's vanishing theorem imply Theorem \ref{main}. 
\end{say}
\begin{say}
The corollaries in Section \ref{sec1} directly follow from 
the main theorem I:~Theorem \ref{main}. We note that 
Corollary \ref{nori} is equivalent to 
the following statement. 
This formulation seems to be more 
useful for various applications. 
\begin{cor}[{cf.~\cite[Corollary 2.5 (iii)]{mustata}}]\label{honyaku} 
Let $X$ be a projective toric variety over $k$ and $L$ an ample 
line bundle on $X$. 
If $D_{j_1}, \cdots, D_{j_r}$ are distinct 
torus invariant prime divisors, then 
$H^i(X, L\otimes \mathcal O_X(-D_{j_1}-\cdots -D_{j_r}))=0$ for every 
$i>0$. 
\end{cor}
\end{say}

Let us go to the proofs 
of the main theorem II:~Theorem \ref{main2}, and Corollary \ref{kavi}. 

\begin{say}[Proof of Theorem \ref{main2}]\label{210}
Let $F:X\to X'$ be the $l$ times multiplication 
map constructed in \ref{21}. Then there 
exist natural split injections 
$\mathcal O_{V'}(\llcorner D'\lrcorner)\to F_*\mathcal O_V(lD)$ and 
$\mathcal O_{V'}(K_{V'}+\ulcorner D'\urcorner)\to F_*\mathcal O_V(K_V+lD)$, 
which are induced by the split injections 
$\mathcal O_{T'}\to F_*\mathcal O_T$ and 
$\Omega^n_{T'}\to F_*\Omega^n_{T}$ (see \ref{24}). 
By pushing them forward to $X'$, we obtain split injections 
$\mathcal O_{X'}(\llcorner D'\lrcorner)\to F_*\mathcal O_X(lD)$ and 
$\mathcal O_{X'}(K_{X'}+\ulcorner D'\urcorner )\to F_*
\mathcal O_X(K_X+lD)$. 
So, we obtain 
\begin{align*}
H^i(X, \mathcal O_{X}(\llcorner D\lrcorner))&\simeq 
H^i(X', \mathcal O_{X'}(\llcorner D'\lrcorner))\\
&\subset 
H^i(X', F_*\mathcal O_X(lD) )\simeq 
H^i(X, \mathcal O_X(lD)) 
\end{align*} 
and 
\begin{align*}
H^i(X, \mathcal O_{X}(K_X+\ulcorner D\urcorner))&\simeq 
H^i(X', \mathcal O_{X'}(K_{X'}+\ulcorner D'\urcorner))\\
&\subset 
H^i(X', F_*\mathcal O_X(K_X+lD) )\\&\simeq 
H^i(X, \mathcal O_X(K_X+lD)).  
\end{align*}  
Thus, Theorem \ref{main2} is obvious. 
\end{say}
\begin{say}[Proof of Corollary \ref{kavi}] 
We take a positive integer $l$ such that $lD$ is integral and Cartier. 
Then $\mathcal O_X(K_X+lD)\simeq \mathcal O_X(K_X)\otimes 
\mathcal O_X(lD)$ since 
$\mathcal O_X(lD)$ is locally free. 
Thus, $H^i(X, \mathcal O_X(K_X)\otimes \mathcal 
O_X(lD))=0$ for 
$i\ne n-\kappa$ (see Theorem \ref{kol} below) and 
$H^i(X, \mathcal O_X(lD))=0$ for 
$i\ne 0$ since $lD$ is a nef Cartier divisor (see, for example, 
\cite[p.74 Corollary]{fulton}). 
This implies the desired vanishing theorems in 
Corollary \ref{kavi}.
\end{say}

\begin{rem}
Note that there are complete toric varieties that have no 
non-trivial nef line bundles (see \cite{f-k} and \cite{fp}). 
\end{rem}

The next remark is due to Nakayama.

\begin{rem}\label{linear}
In Theorem \ref{main2} and Corollary \ref{kavi}, the assumption that $D$ 
is a torus invariant $\mathbb Q$-Weil
divisor on $X$ can be slightly weakened.
It is sufficient to assume that the fractional 
part $\{D\}$ is a torus invariant $\mathbb Q$-Weil
divisor on $X$. We note that 
the integral part $\llcorner D\lrcorner$ is 
always linearly equivalent to a torus invariant Weil divisor on $X$.
Similar modifications work for 
Propositions \ref{vari}, \ref{sp}, Corollary \ref{musumusu}, 
Theorems \ref{43}, and \ref{44} below. 
We leave the details for the readers' exercises.
\end{rem}

\begin{say}
The following proposition is a variant of Theorem \ref{main2}. 

\begin{prop}\label{vari} 
We use the same notation as in {\em{Theorem \ref{main2}}}. 
Let $B$ be a reduced torus invariant Weil divisor on $X$ such that 
$B$ and $\{D\}$ have no common irreducible components. 
If $H^i(X, \mathcal O_X(K_X+B+lD))=0$, 
then $H^i(X, \mathcal O_X(K_X+B+\ulcorner D\urcorner))=0$. 
We further assume that 
$X$ is projective and $D$ is an ample $\mathbb Q$-Cartier $\mathbb Q$-Weil 
divisor. Then 
$H^i(X, \mathcal O_X(K_X+B+\ulcorner D\urcorner))=0$ for 
$i>0$. 
\end{prop} 
The proof is essentially the same as that of Theorem \ref{main2} if we 
use Corollary \ref{nori}. We leave it for the readers' exercise. 
\end{say}

\begin{say}\label{214} 
Let us compare 
Musta\c t\u a's original vanishing 
theorem:~\cite[Theorem 0.1]{mustata} 
with Theorem \ref{main2}. 
The following corollary is nothing but a reformulation 
of Theorem \ref{main2}, which is a slight but 
important generalization 
of Musta\c t\u a's vanishing theorem. 
We note that Doctor Sam Payne independently 
obtained the first part of Corollary \ref{musumusu} 
by another method. 

\begin{cor}\label{musumusu}
Let $X$ be a toric variety defined over $k$ and $D$ a
torus invariant Weil divisor on $X$.
Suppose that we have $E=\sum _{j=1}^{d}a_jD_j$
with $0\leq a_j\leq 1$, where $D_1, \cdots, D_d$ are
distinct torus invariant prime divisors on $X$,
such that $mE$ is an integral Weil divisor for some integer
$m\geq 1$. If $H^i(X, \mathcal O_X(D+m(D+E)))=0$ for
some $i\geq 0$, then $H^i(X, \mathcal O_X(D))=0$. 
Moreover, if $H^i(X, \mathcal O_X(K_X+D+m(D+E)))=0$ 
for some $i\geq 0$, then 
$H^i(X, \mathcal O_X(K_X+D+\ulcorner E\urcorner))=0$. 
\end{cor}

\begin{proof} 
We put $l:=m+1$ and consider a $\mathbb Q$-Weil divisor 
$D^{\dag}:=D+\frac{m}{m+1}E$. 
Then, apply Theorem \ref{main2}. 
We note that $l D^{\dag}=D+m(D+E)$, 
$\llcorner D^{\dag}\lrcorner =D$, 
and $\ulcorner D^{\dag}\urcorner =D+\ulcorner E\urcorner$. 
\end{proof}

\begin{rem}
In Corollary \ref{musumusu}, we do not assume that
$m(D+E)$ is Cartier. So, the first statement is slightly
better than Musta\c t\u a's original
one:~\cite[Theorem 0.1]{mustata}. 
This difference may look very small. 
However, it causes big differences in various 
applications (see Corollary \ref{kavi} and Remark \ref{2200} below). 
The latter statement is new. 
As we saw in Remark \ref{linear}, we do not have to assume 
that $D$ is torus invariant. 
\end{rem}

\begin{rem}\label{2200} 
Let $D$ be a torus invariant $\mathbb Q$-Weil divisor on $X$ such that 
$lD$ is integral. 
If we put $D^{\spadesuit}:=\llcorner D\lrcorner$ 
and $E^{\spadesuit}:=\frac{l}{l-1}\{D\}$, 
and apply Corollary \ref{musumusu} 
to $D^{\spadesuit}$ and $E^{\spadesuit}$ with $m:=l-1$, 
then we can recover Theorem \ref{main2} from 
Corollary \ref{musumusu}. 
To recover Theorem \ref{main2} from 
Musta\c t\u a's theorem:~\cite[Theorem 0.1]{mustata}, 
we have to assume that $m(D^{\spadesuit}+E^{\spadesuit}) 
=lD-\llcorner D\lrcorner$ is Cartier. 
It seems to be a very artificial assumption. 
Thus, I believe that our theorem is much better. 
\end{rem}
\end{say}

\begin{say}
In \cite{vie}, Viehweg obtained his vanishing theorems 
as applications of the Bogomolov type vanishing 
theorem (cf.~\cite[Theorem III]{vie}). 
For toric varieties, we can easily check the following 
Bogomolov type vanishing theorem. 

\begin{thm}[{Bogomolov, $\cdots$}]\label{bogo}  
Let $X$ be a complete toric 
variety defined over a field $k$ and $B$ a reduced torus invariant 
Weil divisor on $X$. 
Let $L$ be a line bundle on $X$ with 
the Iitaka dimension $\kappa (X, L)\geq 0$. 
Then $H^0(X, \widetilde {\Omega}^a_X(\log B)\otimes L^{-1})=0$ 
for any $a\geq 0$ unless $L\simeq \mathcal O_X$. 
\end{thm} 
\begin{proof}
Assume that $H^0(X, \widetilde {\Omega}^a_X(\log B)
\otimes L^{-1})\ne 0$. 
Since $\widetilde {\Omega}^a_{X}(\log B)\subset \wedge^a
M\otimes \mathcal O_X$, we obtain 
$H^0(X, L^{-1})\ne 0$. Therefore, 
$L\simeq \mathcal O_X$ by the 
assumption $\kappa (X, L)\geq 0$. 
\end{proof}

We think that the Kawamata-Viehweg type vanishing theorem for toric 
varieties (cf.~Corollary \ref{kavi}) does not directly follow from 
Theorem \ref{bogo}. 
\end{say}

\begin{say} We close this section with the following three remarks. 

\begin{rem}\label{coxx}
In \cite[Theorem 7.1]{bc}, Corollary \ref{bott} was proved under the 
assumption that the toric variety is $\mathbb Q$-factorial, equivalently, 
has only quotient singularities. 
Batyrev and Cox proved it as a special case of \cite[Theorem 7.2]{bc}. 
We note that we can easily prove \cite[Theorem 7.2]{bc} by 
\cite[Theorem 5.4]{bc}, which is \cite[15.7]{danilov}, and Corollary \ref{bott} 
using induction on $k$ (not on $p-k$). 
For $k$ and $p-k$, see the proof of \cite[Theorem 7.2]{bc}.  
Therefore, we can obtain \cite[Lemma 7.4]{bc} as a corollary 
of the vanishing theorem:~Corollary \ref{bott}. 
Here, we do not pursue this subject anymore since we need the 
Hodge theory. 
\end{rem}

\begin{rem}[Frobenius morphisms]\label{fro}
If $l=p$ and $k$ is a perfect 
field, then $F:V\to V'$ is the relative Frobenius 
morphism and $\psi$ induces 
the inverse Cartier isomorphisms 
$\wedge ^aC^{-1}: \Omega^a_{V'/k}(\log B')\simeq \mathcal H^a(F_*
\Omega^{\bullet}_{V/k}(\log B))$ for any $a\geq 0$. 
All the computations we need were described in \cite[9.14.~Theorem]{ev}. 
We note that 
this technique produces the 
$E_1$-degeneration of the spectral sequence 
$E^{ij}_1=H^j(X, \widetilde \Omega^i_{X}(\log B))\Rightarrow 
\mathbb H^{i+j}(X, \widetilde \Omega^{\bullet}_{X}(\log B))$ 
(see \cite[Remark 1]{btlm}). 
We do not pursue this topic since it was already 
treated in \cite{btlm} and \cite{blickle}. 
\end{rem}

\begin{rem}[Applications of vanishing theorems] 
In Section 4 in \cite{mustata}, Musta\c t\u a 
obtained various results on linear systems on toric 
varieties as applications of 
his vanishing theorem (cf.~\cite[Corollary 2.5 (iii)]{mustata} 
or Corollaries \ref{nori} and \ref{honyaku}). 
In those applications, the considered toric varieties are always 
non-singular. 
In \cite{fujino}, Musta\c t\u a's results in \cite[Section 4]{mustata} 
were reproved and some of them were 
generalized for singular toric varieties. 
See \cite[Section 4 and Remark 3.3]{fujino}. 
However, the proofs in \cite{fujino} are 
quite different from Musta\c t\u a's. 
They depend on the toric Mori theory. 
Note that the foundation of the toric Mori theory was constructed 
without using vanishing theorems (see \cite{reid}, \cite{fs}, \cite{f-3}, and 
\cite{sato}). 
See also \cite[\S 4.~Applications]{sato} for some generalizations of 
Musta\c t\u a's results for the relative setting. 
\end{rem}
\end{say}

\section{Variants of the main vanishing theorems}\label{varivari} 

In this section, we treat slight 
generalizations of the main vanishing theorems. 
We need no new arguments. 

\begin{say}\label{300}
The following theorem is a small generalization of Theorem \ref{main}.
It may be useful in the future. So, we state it here. 

\begin{prop}\label{newnew}
Let $X$ and $B$ be the same as in 
{\em{Theorem \ref{main}}}.
Let $D$ be a $($not necessarily torus invariant$)$ Weil divisor on $X$ 
and $\mathcal E$ a reflexive sheaf on $X$. 
We consider the $l$ times multiplication 
map $F:X\to X'\simeq X$, which was defined in {\em{\ref{21}}}, 
for some positive integer $l$.  
If $H^i(X, (\widetilde \Omega^a_X(\log B)
\otimes F^*\mathcal E\otimes \mathcal O_X(lD))^{**})=0$,
then $H^i(X, (\widetilde \Omega^a_X(\log B)\otimes 
\mathcal E\otimes \mathcal O_X(D))^{**})=0$. 
In particular, 
$H^i(X, (\widetilde \Omega^a_X(\log B)
\otimes \mathcal O_X(lD))^{**})=0$ implies 
$H^i(X, (\widetilde \Omega^a_X(\log B)
\otimes \mathcal O_X(D))^{**})=0$. 
\end{prop} 

We will prove this proposition 
after the proof of Proposition \ref{sp}. 

\begin{rem}
In Proposition \ref{newnew}, 
if $\mathcal E$ is locally free 
and $D$ (resp.~$lD$) is Cartier,  
or $X$ is non-singular, then 
$(\widetilde \Omega^a_X(\log B)
\otimes \mathcal E\otimes \mathcal O_X(D))^{**}
\simeq \widetilde \Omega^a_X(\log B)
\otimes \mathcal E\otimes \mathcal O_X(D)
$ 
(resp.~$(\widetilde \Omega^a_X(\log B)
\otimes F^*\mathcal E\otimes \mathcal O_X(lD))^{**}
\simeq \widetilde \Omega^a_X(\log B)
\otimes F^*\mathcal E\otimes \mathcal O_X(lD)
$) 
in Proposition \ref{newnew}. 
See Remark \ref{2111} (ii) below. 
\end{rem}
\end{say}

\begin{say} 
We treat a similar variant of Theorem \ref{main2} here. 
Doctor Sam Payne independently obtained 
a special case of the 
following theorem under the 
extra assumption that $\mathcal E$ is equivariant. 
I was inspired by his private notes. 

\begin{prop}\label{sp}
We use the same notation as in {\em{Theorem \ref{main2}}}. 
Let $\mathcal E$ be a reflexive sheaf on $X$. 
Let $F:X\to X'\simeq X$ be the 
$l$ times multiplication map as in {\em{\ref{21}}}. 
If $H^i(X, (F^*\mathcal E\otimes \mathcal O_X(lD))^{**})=0$ 
$($resp.~$H^i(X, (F^*\mathcal E
\otimes \mathcal O_X(K_X+lD))^{**})
=0$$)$, then we have $H^i(X, (\mathcal E\otimes 
\mathcal O_X(\llcorner D\lrcorner ))^{**})=0$ 
$($resp.~$H^i(X, (\mathcal E\otimes 
\mathcal O_X(K_X+\ulcorner 
D\urcorner))^{**})=0$$)$. 
\end{prop}

\begin{rem}
If $\mathcal E$ is a locally free sheaf 
or $X$ is non-singular, then we 
do not need to take double duals in Proposition \ref{sp}. 
See Remark \ref{2111} (ii) below. 
If $\mathcal E\simeq \mathcal O_X$, then Proposition 
\ref{sp} is nothing but Theorem \ref{main2}. 
\end{rem}

Before we go to the proofs, we make some remarks on 
reflexive sheaves. 

\begin{rem}\label{2111}
(i) Let $\mathcal F$ be a coherent sheaf on a normal variety $X$. 
Then $\mathcal F^{**}$
denotes the double dual of $\mathcal F$. 
(ii)  Let $\mathcal F_1$ and $\mathcal F_2$ be reflexive 
sheaves on a normal variety $X$. 
Then $(\mathcal F_1\otimes \mathcal F_2)^{**}
\simeq \mathcal F_1\otimes \mathcal F_2$ if 
one of the $\mathcal F_i$ 
is locally free. 
\end{rem}

\begin{proof}[Proof of {\em{Proposition \ref{sp}}}]
Let $V'$ be the Zariski open set of $X'$ as in \ref{26}. 
We take a Zariski open set 
$W'$ of $V'$ such that 
$\mathcal E'$ is locally free on $W'$ and $\codim _{X'}(X'\setminus 
W')\geq 2$. 
Note that $W'$ is not torus invariant when $W'\ne V'$. 
We put $W=F^{-1}(W')\subset V$. Then we obtain 
the following commutative diagram 
$$
\begin{CD}
W@>>>X\\
@V{F}VV @VV{F}V\\
W'@>>>X' 
\end{CD}
$$ 
as in \ref{27}, where the horizontal arrows are 
natural open immersions. 
We have split injections 
$$
\mathcal E'\otimes \mathcal O_{W'}(\llcorner 
D'\lrcorner)\to \mathcal E'\otimes F_*\mathcal O_W(lD)
\simeq F_*(F^*\mathcal E'\otimes \mathcal O_W(lD))
$$ 
and 
\begin{align*}
\mathcal E'\otimes \mathcal O_{W'}(K_{W'}+\ulcorner 
D'\urcorner )&\to \mathcal E'\otimes 
F_*\mathcal O_W(K_W+lD)\\
&\simeq 
F_*(F^*\mathcal E'\otimes \mathcal O_{W}(K_W+lD))
\end{align*}
on $W'$ by \ref{210} and the projection formula. 
By pushing them forward to $X'$, we 
obtain the desired vanishing theorems by the same arguments 
as in \ref{210}. 
\end{proof}
\end{say}

\begin{proof}[Proof of {\em{Proposition \ref{newnew}}}] 
We note that we can replace $D$ by a linearly equivalent torus 
invariant Weil divisor. So, we assume that $D$ is 
torus invariant. 
By the arguments in \ref{26} and \ref{28},
we can check that there exist split injections
$$\Omega^a_{W'}(\log B')\otimes \mathcal E\otimes 
\mathcal O_{W'}(D')\to
F_*(\Omega ^a_W(\log B)\otimes F^*\mathcal E
\otimes \mathcal O_W(lD))$$ for
any $a\geq 0$, where $W'$ (resp.~$W$) is 
the Zariski open set of $V'$ (resp.~$V$) 
defined in the proof of Proposition \ref{sp}. 
Note that $F^*\mathcal O_{V'}(D')\simeq 
\mathcal O_V(lD)$ since $D$ is Cartier 
on $V$. It is because $V$ is non-singular. 
Therefore, by pushing the above split injections to 
$X'$, 
we have split injections
$(\widetilde \Omega^a_{X'}(\log B')\otimes 
\mathcal E\otimes 
\mathcal O_{X'}(D'))^{**}
\to F_*((\widetilde \Omega ^a_X(\log B)
\otimes F^*\mathcal E\otimes \mathcal O_X(lD))^{**})$ for
all $a\geq 0$ (see \ref{27}).
This obviously implies Proposition \ref{newnew}.
\end{proof} 
By applying Proposition 
\ref{newnew} in place of Theorem \ref{main}, 
some vanishing theorems in this paper can be generalized slightly. 
We leave the details for the readers' exercise. 
We also leave to the interested readers the pleasure of 
combining 
the latter part of Proposition \ref{sp} with 
Proposition \ref{vari}. 

\section{Koll\'ar's injectivity theorem}\label{sec3} 

In this section, we treat Koll\'ar's injectivity theorem 
(cf.~\cite[Theorem 2.2]{koll}) for toric varieties.
It is an application of Corollary \ref{kodaira}. 
 
\begin{thm}\label{kol}
Let $X$ be a complete toric variety defined over $k$ and $L$ a nef line 
bundle on $X$. 
Let $s$ be a non-zero holomorphic section of $L^l$, where 
$l\geq 0$. 
Then 
$$
\times s:H^i(X, \mathcal O_X(K_X)\otimes 
L^m)\to H^i(X, \mathcal O_X(K_X)\otimes L^{m+l}) 
$$  
is injective for any $m\geq 1$ and $i\geq 0$, 
where $\times s$ is the morphism induced by the tensor product with $s$. 
More precisely, $H^i(X, \mathcal O_X(K_X)\otimes L^m)=0$ for any $m\geq 1$ 
when $i\ne n-\kappa$. 
Here, $n=\dim X$ and $\kappa =\kappa (X, L)$. 
\end{thm}

The following lemma is well known. The readers can find it in 
any text book on the toric geometry (see, for example, \cite[p.76 
Proposition and p.89 Proposition]{fulton}).  

\begin{lem}\label{classical}
Let $f:X\to Y$ be a proper birational toric morphism. 
Then $R^if_*\mathcal O_X=0$ and $R^if_*\mathcal O_X(K_X)=0$ for all 
$i>0$, 
$f_*\mathcal O_X\simeq \mathcal O_Y$, and $f_*\mathcal O_X(K_X)\simeq 
\mathcal O_Y(K_Y)$. 
\end{lem} 

The next lemma is a slight generalization of Lemma \ref{classical}. 

\begin{lem}\label{new}
Let $f:X\to Y$ be a proper surjective toric morphism with connected fibers. 
Then $R^if_*\mathcal O_X=0$ for $i>0$ and 
$f_*\mathcal O_X\simeq \mathcal O_Y$. 
Moreover, $R^{n-m}f_*\mathcal O_X(K_X)\simeq 
\mathcal O_Y(K_Y)$ and $R^if_*\mathcal O_X(K_X)=0$ for $i\ne n-m$, 
where $n=\dim X$ and $m=\dim Y$. 
\end{lem} 

\begin{proof}[Sketch of the proof] 
The former statement is an exercise if we use Lemma \ref{classical}. 
For the proof, see, for example, \cite[Theorem 3.2]{ishida}. 
The latter part follows from the Grothendieck duality and the former 
statement. 
\end{proof} 

Lemma \ref{new} implies that Koll\'ar's torsion-freeness (cf.~\cite
[Theorem 2.1 (i)]{koll}) is obvious 
for toric varieties and Koll\'ar's vanishing theorem (cf.~\cite[Theorem 
2.1 (iii)]{koll}) is a special case 
of Corollary \ref{kodaira} in the toric geometry. 

\begin{proof}[{Proof of {\em{Theorem \ref{kol}}}}] 
Since $L$ is nef, there exists a proper surjective toric 
morphism with connected fibers $f:X\to Y$ such that 
$L\simeq f^*H$, where $H$ is an ample line bundle on $Y$. 
By the definition of $\kappa$, we have $\dim Y=\kappa$. 
We consider the spectral sequence $H^i(Y, R^jf_*\mathcal O_X(K_X)
\otimes H^b)\Rightarrow H^{i+j}(X, \mathcal O_X(K_X)\otimes L^b)$ for any 
integer $b$. 
By Lemma \ref{new}, we obtain $H^i(Y, \mathcal O_Y(K_Y)\otimes 
H^b)\simeq H^{i+n-\kappa}(X, \mathcal O_X(K_X)\otimes L^b)$. 
Therefore, we have $H^{n-\kappa}(X, \mathcal O_X(K_X)\otimes L^b)\simeq 
H^0(Y, \mathcal O_Y(K_Y)\otimes H^b)$ and 
$H^{i}(X, \mathcal O_X(K_X)\otimes L^m)=0$ for $i\ne n-\kappa$ 
and $m\geq 1$ by Corollary \ref{kodaira}. 
Note that $H^0(X, L^l)\simeq H^0(Y, H^l)$. 
So, there exists a non-zero $t\in H^0(Y, H^l)$ such that $s=f^*t$. 
Thus,  
$
\times s:H^{n-\kappa}(X, \mathcal O_X(K_X)\otimes 
L^m)\to H^{n-\kappa}(X, \mathcal O_X(K_X)\otimes L^{m+l})$ 
is nothing but 
$
\times t :H^0(Y, \mathcal O_Y(K_Y)\otimes 
H^m)\to H^0(Y, \mathcal O_Y(K_Y)\otimes H^{m+l}) 
$. 
Therefore, $\times s$ is injective since $\times t$ is injective. 
\end{proof}

\begin{say}
As we saw in Theorem \ref{kol}, 
the Kodaira type vanishing theorem (cf.~Corollary \ref{kodaira})
holds for nef and big line bundles. 
However, the Norimatsu type vanishing theorem (cf.~Corollary \ref{nori})
does not always hold for nef and big line bundles by the next example.

\begin{ex}
In this example, we assume $k=\mathbb C$, 
the complex number field, for simplicity.
Let $P\in \mathbb P^2$ be a torus invariant 
closed point 
and let $f:X\to \mathbb P^2$ be the
blow-up at $P$. Let $B$ be the $f$-exceptional curve on $X$. 
Then we obtain
$0\to \mathcal O_X(K_X)\to \mathcal 
O_X(K_X+B)\to \mathcal O_B(K_B)\to 0$ by adjunction.
By applying $R^if_*$, we obtain 
$f_*\mathcal O_X(K_X+B)\simeq 
\mathcal O_{\mathbb P^2}(K_{\mathbb P^2})$
and $R^1f_*\mathcal O_X(K_X+B)\simeq 
\mathbb C(P)$ since $R^if_*\mathcal O_X(K_X)=0$ for
$i>0$. We put $H=\mathcal O_{\mathbb P^2}(1)$ 
and $L=f^*H$. Note that
$L$ is nef and big. Then, 
by the Leray spectral sequence,
we have the following exact sequence: 
$0\to H^1(\mathbb P^2, f_*\mathcal O_X(K_X+B)\otimes H)\to
H^1(X, \mathcal O_X(K_X+B)\otimes L)\to
H^0(\mathbb P^2, R^1f_*\mathcal O_X(K_X+B)\otimes H)\to
H^2(\mathbb P^2, f_*\mathcal O_X(K_X+B)\otimes H)\to \cdots .$
Since the first and the last terms are zero, $H^1(X, \mathcal O_X(K_X+B)\otimes L)
\simeq H^0(\mathbb P^2, \mathbb C(P)\otimes H)\simeq \mathbb C$. 
\end{ex}
\end{say}

\section{Appendix:~Relative vanishing theorems}\label{sec4}
In this appendix, we state relative vanishing theorems explicitly 
for future uses. 
All the vanishing theorems easily follow from the main theorems 
and their proofs. We only give a proof of Theorem \ref{43} for the 
readers' convenience. The others are similar and easier to prove. 

\begin{say}
Let $f:X\to Y$ be a proper 
surjective toric morphism and $B$ a reduced 
torus invariant Weil divisor on $X$. We put $n=\dim X$ 
and $m=\dim Y$. 

\begin{thm}\label{th-a}
Let $L$ be an $f$-ample line bundle 
on $X$. 
Then we have $R^if_*(\widetilde \Omega^a_{X}(\log B)\otimes 
L)=0$ for $i>0$. 
In particular, $R^if_*(\widetilde \Omega^a_X\otimes L)=0$, 
$R^if_*(\mathcal O_X(K_X+B)\otimes L)=0$, and 
$R^if_*(\mathcal O_X(K_X)\otimes L)=0$ for $i>0$. 
\end{thm} 

\begin{thm}\label{43}
Let $D$ be a torus invariant $\mathbb Q$-Weil divisor on $X$. 
Assume that $D$ is $\mathbb Q$-Cartier and $f$-nef 
with the relative Iitaka dimension 
$\kappa (X/Y, D)=\kappa$. 
Then $R^if_*\mathcal O_X(K_X+\ulcorner D\urcorner )=0$ for 
$i\ne n-m-
\kappa$ and $R^if_*\mathcal O_X(\llcorner D\lrcorner )=0$ for 
$i\ne 0$. 
\end{thm} 

\begin{proof}
Let $l$ be a positive integer such that 
$lD$ is integral and Cartier. 
We note that $f^*f_*\mathcal O_X(lD)\to \mathcal O_X(lD)$ is 
surjective since $lD$ is an $f$-nef 
Cartier divisor on $X$ (see, for example, \cite[Chapter 
IV, 1.13 Lemma]{nakayamabon}). 
We can assume that $Y$ is affine since the problem is 
local. 
By \ref{210}, 
it is sufficient to prove that 
$R^if_*(\mathcal O_X(K_X)\otimes 
\mathcal O_X(lD))=0$ for $i\ne n-m-\kappa$ and 
$R^if_*\mathcal O_X(lD)=0$ for $i\ne 0$. 
First, we prove $R^if_*\mathcal O_X(lD)=0$ for 
$i>0$. 
In this case, $R^if_*\mathcal O_X(lD)\simeq 
H^i(X, \mathcal O_X(lD))=0$ by \cite[p.74 Corollary]{fulton} 
since $\mathcal O_X(lD)$ is generated by its global 
sections and the support of 
the fan associated to 
$X$ is convex. 
Next, we prove 
$R^if_*(\mathcal O_X(K_X)\otimes 
\mathcal O_X(lD))=0$ for $i\ne n-m-\kappa$. 
Let $g:X\to Z$ be a proper surjective 
toric morphism  over $Y$ with connected fibers such 
that $\mathcal O_X(lD)\simeq g^*\mathcal O_Z(H)$, 
where $H$ is a Cartier divisor on $Z$ which is ample over $Y$. 
We note that $\dim Z=m+\kappa$. 
By Lemma \ref{new} and 
Leray's spectral sequence, we obtain 
$R^if_*(\mathcal O_X(K_X)\otimes \mathcal O_X(lD))\simeq 
R^{i-(n-m-\kappa)}h_*(\mathcal O_Z(K_Z)\otimes 
\mathcal O_Z(H))$, 
where $h:Z\to Y$. 
In particular, $R^if_*(\mathcal O_X(K_X)\otimes \mathcal O_X(lD))
=0$ for $i<n-m-\kappa$. 
By the same arguments as in \ref{28}, we have 
$R^{i-(n-m-\kappa)}h_*(\mathcal O_Z(K_Z)\otimes 
\mathcal O_Z(H))\subset 
R^{i-(n-m-\kappa)}h_*(\mathcal O_Z(K_Z)\otimes 
\mathcal O_Z(l'H))$ for any positive integer $l'$. 
If $l'\gg 0$, then $R^{i-(n-m-\kappa)}h_*(\mathcal 
O_Z(K_Z)\otimes \mathcal O_Z(l'H))=0$ for 
$i>n-m-\kappa$ by Serre's vanishing theorem. 
Thus, $R^{i-(n-m-\kappa)}h_*(\mathcal O_Z(K_Z)
\otimes \mathcal O_Z(H))=0$ for 
$i>n-m-\kappa$. 
Therefore, we obtain the desired vanishing theorems. 
\end{proof}

\begin{thm}\label{44}
Let $D$ be an $f$-ample 
$\mathbb Q$-Cartier torus invariant $\mathbb Q$-Weil divisor on $X$ such that 
$B$ and $\{D\}$ have no common irreducible components. 
Then $R^if_*\mathcal O_X(K_X+B+\ulcorner D\urcorner )=0$ for $i>0$. 
\end{thm} 

We note that Theorem \ref{kol} can be generalized for the relative 
setting if we use Theorem \ref{th-a} instead of Corollary \ref{kodaira}. 
\end{say}
\ifx\undefined\bysame
\newcommand{\bysame|{leavemode\hbox to3em{\hrulefill}\,}
\fi

\end{document}